\documentclass[11pt]{amsart}
\usepackage{graphicx}
\usepackage{float}
\usepackage{amssymb,amsmath,amsthm}
\usepackage{epstopdf}
\usepackage{color}
\usepackage{mdframed}
\usepackage{empheq}
\usepackage{bm}
\usepackage{hyperref}
\usepackage{mathabx} % for \widecheck
\usepackage{enumitem}
\usepackage{xcolor}
\usepackage[linesnumbered,ruled,vlined]{algorithm2e}
\usepackage{subcaption}
\usepackage{multirow}

\makeatletter
\renewcommand\l@paragraph[2]{}
\renewcommand\l@subparagraph[2]{}
\makeatother
\title[Cost Attribution and Risk-Averse]{Cost Attribution And Risk-Averse Unit Commitment In Power Grids Using Integrated Gradient}
\author{Xinshuo Yang, Ren\'e Carmona, Ronnie Sircar}
\address{Department of Operations Research \& Financial Engineering\\
Princeton University\\
Princeton, NJ 08544, USA}
\date{\today}                   

\begin{document}
\hfill

\begin{abstract}

This paper introduces a novel approach to addressing uncertainty and associated risks in power system management, focusing on the discrepancies between forecasted and actual values of load demand and renewable power generation. By employing Economic Dispatch (ED) with both day-ahead forecasts and actual values, we derive two distinct system costs, revealing the financial risks stemming from uncertainty. We present a numerical algorithm inspired by the Integrated Gradients (IG) method to attribute the contribution of stochastic components to the difference in system costs. This method, originally developed for machine learning, facilitates the understanding of individual input features' impact on the model's output prediction. By assigning numeric values to represent the influence of variability on operational costs, our method provides actionable insights for grid management. As an application, we propose a risk-averse unit commitment framework, leveraging our cost attribution algorithm to adjust the capacity of renewable generators, thus mitigating system risk. Simulation results on the RTS-GMLC grid demonstrate the efficacy of our approach in improving grid reliability and reducing operational costs. 

\end{abstract}

\maketitle


%%%%%%%%%%%%%%%%%%%%%%%%%%%%%
\section{\textbf{Introduction}}
\label{se:introduction}
%%%%%%%%%%%%%%%%%%%%%%%%%%%%

In the pursuit of a more sustainable and environmentally responsible energy landscape, the world is experiencing a remarkable shift towards high levels of renewable energy penetration. The relentless expansion of wind, solar, and other clean energy technologies has brought us closer to reducing our carbon footprint and combating climate change. However, the integration of renewable sources introduces an unprecedented level of complexity into many fundamental aspects of power system management. Unlike their conventional counterparts, renewable resources are inherently variable, contingent upon factors like weather patterns, sunlight availability, and wind speed. Consequently, new sources of randomness as well as inherent risk factors should be introduced into Unit Commitment (UC) and Economic Dispatch (ED), two crucial components of modern power system management. In the short term operational scheduling of electric power systems, UC involves the strategic selection of power generation units, determining which sources, whether they are conventional or renewable, should be dispatched to meet electricity demand over a specific time horizon. ED, on the other hand, optimizes the allocation of generation output among these chosen units in order to minimize costs while ensuring the reliable provision of electricity.

\vskip 4pt

Traditionally, UC models take a deterministic approach in which the parameters and variables are assumed to be known exactly \cite{padhy2004unit}. Due to their deterministic nature, these models are not capable of capturing variability and uncertainty. In recent years, there has been significant research in modeling and addressing the reliability of UC under uncertainty. Among those works, stochastic programming and robust optimization are two common methodologies frequently used. 
In stochastic programming models, the expected operating cost is minimized under a set of representative scenarios to meet the demand while satisfying the operational constraints \cite{tuohy2009unit, wu2007stochastic, zheng2014stochastic, spyrou2024assess}. These representative scenarios may sometimes encompass highly rare events, leading to solutions that come with substantial costs. To tackle this issue, risk-averse UC models were developed with additional constraints to manage the risk exposures through the use of risk measures, conditional value at risk (CVaR) being the most widely used \cite{jabr2005robust, huang2014two, zhang2018conditional, asensio2015stochastic}. Despite this, stochastic UC models can present computational difficulties, primarily arising from their inherent complexity, particularly when a substantial number of scenarios are integrated into the model.

\vskip 4pt

Unlike models based on stochastic programming, robust UC models aim to address uncertainty by considering only the range of uncertainty, rather than relying on information about underlying probability distributions. Instead of minimizing the total expected cost, as is done in stochastic UC, robust UC focuses on minimizing the worst-case cost across all potential outcomes of uncertain parameters \cite{bertsimas2012adaptive, jiang2011robust, wang2013two}. While robust UC solutions tend to be quite conservative, they have the advantage of avoiding the computational complexity associated with incorporating a large number of Monte Carlo scenarios.

\vskip 4pt

In the present paper, we address the aforementioned uncertainty and its associated risks from a different perspective. The primary source of uncertainty stems from the deviations between forecasted and actual values, whether we are examining load demand or renewable power generation. By running ED with two sets of parameters—the day-ahead forecasts and the actual values—we obtain two distinct system costs. The difference between these costs represents the financial risks associated to these uncertainties. We present a numerical algorithm inspired by a widely recognized method known as the Integrated Gradients (IG) \cite{sundararajan2017axiomatic} to assess the individual contributions to each stochastic component to the difference of system costs. Originally introduced as an attribution mechanism in machine learning, the IG is typically used to attribute the significance or
contribution of individual input features (e.g., pixels in an image or words in text) on a model’s output
prediction, especially in deep neural network learning. The IG chooses a baseline input to ideally represent the absence or neutrality of input features. It then calculates the contribution of each feature to the model's output by integrating the gradients of the model's prediction with respect to the input along a straight path from a baseline to the target input. In the proposed algorithm, the power production cost (ED) model is treated as a function that take input consisting of the random quantities such as load demand and renewable generations, and yields the system cost.  We designate the day-ahead forecast as the baseline input, while the actual values serve as the target input. The algorithm will assign a numeric monetary value to represent the influence of this variability on the grid operational costs. 

Certainly, these values can be used as mere informational tools in various contexts. As an example of application of the cost attribution algorithm, we propose a risk averse unit commitment framework to adjust the capacity of renewable generators using a large number of Monte Carlo scenarios as the target inputs in the attribution method. This
framework considers the distribution of the cost attribution of individual generator through the Monte Carlo scenarios. It
results in a reduction of their generation capacity (maximum generation level) if a certain portion of
that capacity is deemed unlikely to be delivered and could potentially lead to elevated operational costs. By integrating this adjustment mechanism into the Unit Commitment (UC) optimization process, we can avoid relying on highly uncertain renewable productions without sacrificing the benefits of inexpensive and eco-friendly renewable energy, thereby reducing the overall risk in grid operations.

\vskip 4pt

This work is part of research from the Princeton team ORFEUS\footnote{ORFEUS: Operational Risk Financialization
of Electricity Under Stochasticity, \url{orfeus.princeton.edu}}, funded by ARPA-E under its program PERFORM\footnote{PERFORM—Performance-based Energy Resource Feedback, Optimization, and Risk Management}. A related paper that uses Shapley values for renewables reliability risk allocation is \cite{riskalloc}.

The rest of the paper is organized as follows: in Section \ref{se:risk_alloc} we provide a comprehensive description of our problem setup and model configuration and introduce a numerical algorithm designed for attributing operational costs to individual components within this setup. We outline certain desirable properties of the attribution algorithm along with implementation techniques. Section \ref{se:case_study} is devoted to a case study for the illustrative RTS-GMLC grid. We report simulation results and performance of our cost attribution algorithm. In Section \ref{se:risk_averse_uc}, we describe a framework on how to adjust the generation capacity of renewable generators using Monte Carlo scenarios and present yearly simulation results on the RTS-GMLC grid, highlighting its efficacy in mitigating system risk concerning  loss of load.  We conclude with a short recall of the main contributions of the paper in Section \ref{se:conclusion}.

%%%%%%%%%%%%%%%%%%%%%%%%%%%%%
\section{\textbf{Cost Attribution}}
\label{se:risk_alloc}

\subsection{Integrated Gradient}

Integrated Gradients (IG) is a technique widely used to attribute the predictions of a machine learning model to its input features. It is often used to interpret and explain the predictions of complex models, especially deep neural networks. The main goal of IG is to understand the contribution of each feature in the input to the final prediction made by the model. The idea is to compute the integral of the partial derivatives of the model's output with respect to the input features while gradually transitioning from a baseline input (usually a neutral or zero input) to the actual input for which one wants to explain the prediction. This process helps in attributing the prediction incrementally to each feature while considering their contributions in the context of the overall prediction.
This attribution method is reminiscent of the way sensitivities of derivative prices and portfolio values are understood and computed in financial engineering.

Mathematically speaking, the premise is the first order Taylor expansion, also known as the fundamental theorem of calculus. Suppose we have a function $F:\mathbb{R}^n\rightarrow\left[0,1\right]$ that represents a deep neural network, an arbitrary input vector $x\in\mathbb{R}^n$ and a baseline input vector $x'\in\mathbb{R}^n$. The IG provide a quantitative measure of the contribution of each input feature $x_i$ to the model's prediction $F\left(x\right)$. Higher attributions indicate greater influence on the prediction, while lower attributions suggest lesser impact.
The attribution for the $i$-th feature is obtained by integrating the gradients along the straight path from the baseline input $x'$ to the actual input $x$:

\begin{equation}
    I_i\left(x\right) \coloneqq \left(x_i - x'_i\right)\int_{\lambda=0}^1\frac{\partial F\left(x' + \lambda \left(x - x'\right) \right)}{\partial x_i}d\lambda.
\end{equation}

 This integration is typically performed using numerical methods like the trapezoidal rule or Simpson's rule.  One important property of Integrated Gradients is that the sum of the attribution scores across all features equals the difference in the model's output between the actual input and the baseline, i.e.

\begin{equation}
    \sum_{i=1}^nI_i\left(x\right)=F\left(x\right) - F\left(x'\right).    
\end{equation}

\subsection{Production Cost Model}
We adapt the technique of  IG to the problem of cost attribution in power production cost models (PCM) using a single scenario of electricity demand and renewable production. Production cost modeling is a technique to simulate and analyze the economic and operational characteristics of electricity generation and distribution. It involves creating mathematical models that represent the various components of an electrical power system, including power plants, transmission lines, distribution networks, and demand profiles. Modern production cost modeling involves simulations of electricity markets based on forms of Unit Commitment (UC) and Economic Dispatch (ED) for generation and load scheduling. While the mathematical formulations in UC and ED can vary widely depending on the software package and its options \cite{knueven2020mixed}, the cost attribution method proposed here is agnostic to the choice of their formulation.  Power cost modeling  typically  requires alternating day-ahead UC and real-time ED market optimizations over a fixed horizon, whether it is a day, month or year. Our method can be used for most of the common simulation frameworks, however, for the purpose of illustration, we consider a day-ahead and real-time optimization cycle framework executed on a daily basis. The day-ahead simulation consists of a single UC optimization step that determines $48$ hours of system operations subject to load and renewable production forecasts. Subsequently the real-time simulation consists of $24$ ED on the hourly timescale, which re-optimize system operations subject to actual load and renewable power productions at each hour with a $h$-hour look-ahead horizon. A common choice of $h$ can vary from $1$ to $4$ hours. Because many baseload generators are inflexible at short timescales, the real-time simulation treats day-ahead baseload UC profiles as fixed constraints. 

Given the actual and day-ahead forecasts of the load demand and renewable production, our approach to the operational cost attribution is top-down. We first solve the day-ahead UC using the load and renewable generation forecasts. Following the operational constraints in the UC profile, we then assess the system-level costs by solving 24 hourly EDs twice, once with the actual data and the second time with the day-ahead forecasts as the actual data. Finally, we attribute the system-wide cost difference between the two ED to individual load and production assets. 

To be more specific, we consider a system with $M$ dispatchable (conventional) power plants and denote their power generation at time $t$ by $\mathbf{p}_t=\left(p_t^1,p_t^2,\cdots,p_t^M\right)^T\in\mathbb{R}_{+}^M$. For each $\tau=1,2,\cdots,24$, the ED model with  look-ahead horizon $h$ can be written as, 
\begin{equation}
\label{eq:sced}
\begin{aligned}
\min \quad & \sum_{t=\tau}^{\tau+h} \left( \sum_{g=1}^M C_{g}^t\left(p_t^g\right) + C_{p}^t\left(\mathbf{p}_t\right)\right)\\
\textrm{s.t.} \quad & \left(\mathbf{p}_{\tau}^T, \mathbf{p}_{\tau+1}^T, \cdots, \mathbf{p}_{\tau+h}^T\right)\in\Omega   \\
\end{aligned}
\end{equation}
where $C_g^t$ is the operational cost associated with generation unit $g$ at time $t$, $C_p^t$ is a penalty term accounting for loss of load, unmet reserve requirement, etc. at time $t$, and $\Omega$ is the feasible set representing the system constraints (load balance, generation limits, ramping capacities etc). The costs $C_g^t$ and $C_p^t$ as well as the feasible set $\Omega$ are determined by the grid configuration, and so are the load demand profile, renewable generation data and the initial grid operating condition. In other words, given the load demand profile, renewable generation data and initial state of the dispatchable units $\mathbf{p}_0$, one can solve (\ref{eq:sced}) and obtain an optimal operational cost (assuming that the model is feasible). 

Without loss of generality, we set the look-ahead horizon to $h=1$, and we assume that the system has $L$ load buses, and $N$ renewable power plants with load demand profile and power generation capacity at time $t$ given by $\mathbf{d}_t=\left(d_t^1,d_t^2,\cdots,d_t^L\right)^T\in\mathbb{R}_{+}^L$ and $\mathbf{q}_t=\left(q_t^1,q_t^2,\cdots,q_t^N\right)^T\in\mathbb{R}_{+}^N$ respectively. Now the ED model at time $\tau$ can be viewed as a mapping $F_{\tau}$:

\begin{equation}
 \mathbb{R}_{+}^M\times\mathbb{R}_{+}^L\times\mathbb{R}_{+}^N\ni
 (\mathbf{p}_{\tau-1}, \mathbf{d}_{\tau},
        \mathbf{q}_{\tau})\rightarrow
     F_{\tau}\begin{pmatrix}
        \mathbf{p}_{\tau-1}\\ \mathbf{d}_{\tau}\\
        \mathbf{q}_{\tau}
    \end{pmatrix}\in\mathbb{R}_{+},
\end{equation}
providing the overall cost to the system. The goal of the proposed method is to quantify accurately the contributions of each individual load and each generation asset to the system-wide cost difference computed with the day-ahead forecasts and the actual data, namely, for $\tau=1,2,\cdots,24$.
\begin{equation}
    C_{\tau}^{act} - C_{\tau}^{fcst} := F_{\tau}\begin{pmatrix}
        \mathbf{p}_{\tau-1}^{act}\\
        \mathbf{d}_{\tau}^{act}\\
        \mathbf{q}_{\tau}^{act}
    \end{pmatrix} - F_{\tau}\begin{pmatrix}
        \mathbf{p}_{\tau-1}^{fcst}\\
        \mathbf{d}_{\tau}^{fcst}\\
        \mathbf{q}_{\tau}^{fcst}
    \end{pmatrix}.
\end{equation}

\subsection{Cost Attribution of PCM} We assume that $F_{\tau}$ is differentiable almost everywhere, and let $\mathcal{C}$ be a straight line given by 
\begin{equation}
\vec{r}_{\tau}\left(\lambda\right) = \begin{pmatrix}
            \mathbf{p}_{\tau-1}^{fcst} + \lambda\left(\mathbf{p}_{\tau-1}^{act} - \mathbf{p}_{\tau-1}^{fcst}\right)\\
            \mathbf{d}_{\tau}^{fcst} + \lambda\left(\mathbf{d}_{\tau}^{act} - \mathbf{d}_{\tau}^{fcst}\right)\\
            \mathbf{q}_{\tau}^{fcst} + \lambda\left(\mathbf{q}_{\tau}^{act} - \mathbf{q}_{\tau}^{fcst}\right)
            \end{pmatrix}, \quad 0\leq\lambda\leq1
\end{equation}
then by the fundamental theorem of calculus for line integral, we can write

\begin{align*} 
 F_{\tau}\begin{pmatrix}
        \mathbf{p}_{\tau-1}^{act}\\
        \mathbf{d}_{\tau}^{act}\\
        \mathbf{q}_{\tau}^{act}
    \end{pmatrix} - F_{\tau}\begin{pmatrix}
        \mathbf{p}_{\tau-1}^{fcst}\\
        \mathbf{d}_{\tau}^{fcst}\\
        \mathbf{q}_{\tau}^{fcst}
    \end{pmatrix}& =  \int_{\mathcal{C}}\nabla F_{\tau}
     \begin{pmatrix}
         \mathbf{p}_{\tau-1}\\
         \mathbf{d}_{\tau}\\
         \mathbf{q}_{\tau}
     \end{pmatrix}\cdot d\vec{r} \\ 
 = &  \begin{pmatrix}
            \mathbf{p}_{\tau-1}^{act} - \mathbf{p}_{\tau-1}^{fcst}\\
            \mathbf{d}_{\tau}^{act} - \mathbf{d}_{\tau}^{fcst}\\
            \mathbf{q}_{\tau}^{act} - \mathbf{q}_{\tau}^{fcst}
            \end{pmatrix} \cdot \int_{0}^1 \nabla F_{\tau}\left(\vec{r}_{\tau}\left(\lambda\right)\right)d\lambda.
\end{align*}

In particular, the IG assigns to the initial state of dispatchable generator $m$ the attribution:
\begin{equation}
    \label{eq:attri init}
    C_{init}^{m,\tau} = \left(p_{\tau-1}^{m,act}-p_{\tau-1}^{m,fcst}\right)
    \int_{0}^1\frac{\partial F_{\tau}}{\partial p_{\tau-1}^m}
    \left(\vec{r}_{\tau}\left(\lambda\right)\right)d\lambda
\end{equation}
to the load asset $l$,
\begin{equation}
    \label{eq:attri load}
    C_{load}^{l,\tau} = \left(d_{\tau}^{l,act}-d_{\tau}^{l,fcst}\right)
    \int_{0}^1\frac{\partial F_{\tau}}{\partial p_{\tau}^l}
    \left(\vec{r}_{\tau}\left(\lambda\right)\right)
            d\lambda
\end{equation}
and to the renewable generator $n$
\begin{equation}
    \label{eq:attri renew}
    C_{renew}^{n,\tau} = \left(q_{\tau}^{n,act}-q_{\tau}^{n,fcst}\right)
    \int_{0}^1\frac{\partial F_{\tau}}{\partial q_{\tau}^n}
    \left(\vec{r}_{\tau}\left(\lambda\right)\right)
            d\lambda.
\end{equation}
We also have
\begin{equation}
    \label{eq:complete axiom}
    C_{act}^{\tau} - C_{fcst}^{\tau} = \sum_{m=1}^M C_{init}^{m,\tau} + \sum_{l=1}^L C_{load}^{l,\tau} + \sum_{n=1}^N C_{renew}^{n,\tau}.
\end{equation}

Equation (\ref{eq:complete axiom}) is called the {\itshape{Completeness Axiom}}. Basically, it says that the attributions to the initial states of the dispatchable generators, the load demands and the renewable generations add up to the difference between the cost of operating the grid under the day-ahead forecasts and the actual data. The assumption of differentiability of $F_{\tau}$ is obviously a strong condition and can be difficult to check in practice, especially as it depends upon the grid configuration and the ED model formulation. Nevertheless, we observe empirically that indeed the Completeness Axiom holds in our experiments (see Table \ref{ta:alloc_summary}). 

Another important property which is desirable for attribution methods is {\itshape{Symmetry Preserving}}. In the present context of production cost modeling this property can be understood in the following way. A function is said to be symmetric with respect to two input variables if the value of the function remains unchanged when the values of the variables are interchanged. Accordingly, an
attribution method is said to be {\itshape{Symmetry Preserving}}, if whenever the function is symmetric with respect to two input variables, these input variables receive identical attributions. For a power grid system, {\itshape{Symmetry Preserving}} suggests that if two assets have the same impact in terms of the operational cost under all circumstances, they should receive identical attributions. For instance, if two wind farms are connected to the same bus and always produce the same amounts of power, their cost attributions by the proposed method should be identical. At an intuitive level, such a requirement is desirable for a practical attribution to be fair and equitable, and it holds in the attribution method introduced in this paper.

\vskip 2pt
The original introductions and implementations of the IG also mentioned a few other desirable properties, such as \emph{Sensitivity, Linearity, Implementation Invariance, etc}, however, they are not relevant to models of grid management, and we chose to ignore them in our context.

\subsection{Computation of the Gradient}

The integrals in (\ref{eq:attri init})-(\ref{eq:attri renew}) can be efficiently approximated by adaptive quadrature, e.g. adaptive trapezoidal rule. Adaptivity  usually requires less than $100$ function (integrand) evaluations to reach a desire error threshold (within $5\%$). 

For the gradient, noting that the ED model we consider is a linear program, we compute the partial derivative w.r.t. a quantity of interest by the dual value of the constraint corresponding to that quantity. For example, let us assume that at time $\tau$ the generation capacity of the $n$-th renewable generator is $\bar{q}_{\tau}^{n}$.
To compute the partial derivative $\frac{\partial F_{\tau}}{\partial q_{\tau}^n}$, we first add a constraint 
\begin{equation}
    \label{eq:renew constr}
    q_{\tau}^n \leq \bar{q}_{\tau}^{n}
\end{equation}
to the ED model (\ref{eq:sced}). Once the ED model is optimized, we obtain the desired partial derivative as the dual value of the constraint (\ref{eq:renew constr}).This is what is usually called the \emph{shadow price}

%%%%%%%%%%%%%%%%%%%%%%%
\section{\textbf{Case Study: the RTS-GMLC Grid Model}}
\label{se:case_study}
In this section, we present numerical results of the cost attribution method for the RTS-GMLC grid model whose characteristics can be found at \url{https://github.com/GridMod/RTS-GMLC}. The RTS-GMLC grid contains $73$ buses, $157$ generators (including $4$ wind and $57$ solar stations) and $120$ transmission lines. We use the production cost modeling tool Vatic v0.4.1-a1 (\url{https://github.com/PrincetonUniversity/Vatic/releases/tag/v0.4.1-a1}) and GUROBI solver v10.0.1 \cite{gurobi} to run UC and ED optimization. Vatic applies mixed-integer linear programming optimization to power grid formulations created using the software package EGRET \cite{egret}. Our Python implementation of the cost attribution algorithm with detailed examples is available at \url{https://github.com/PrincetonUniversity/PGrisk}\cite{Yang_PGrisk_2023}. For all examples presented in this section, we choose a maximum of $2^{12}=4096$ for the adaptive trapezoidal quadrature nodes with $5\%$ error threshold. Below, we summarize the formulations in our experiments in Table \ref{ta:formulation} and refer to \cite{egret, knueven2020mixed} for more details.

\begin{table}[hbt!]
\begin{tabular}{ c|c|c } 
\label{ta:formulation}
  & UC & ED \\ 
  \hline
 Status & \texttt{garver\_3bin\_vars} & \texttt{garver\_3bin\_vars} \\ 
 \hline
 Power & \texttt{garver\_power\_vars} & \texttt{garver\_power\_vars} \\ 
 \hline
 Reserve & \texttt{garver\_power\_avail\_vars} & \texttt{MLR\_reserve\_vars} \\ 
 \hline
 Generation & \texttt{pan\_guan\_gentile\_KOW\_generation\_limits} & \texttt{MLR\_generation\_limits} \\ 
 \hline
 Ramping & \texttt{damcikurt\_ramping} & \texttt{damcikurt\_ramping} \\ 
 \hline
 Production & \texttt{KOW\_production\_costs\_tightened} & \texttt{CA\_production\_costs} \\ 
 \hline
 Up/Down & \texttt{rajan\_takriti\_UT\_DT} & \texttt{rajan\_takriti\_UT\_DT} \\ 
 \hline
 StartUp & \texttt{KOW\_startup\_costs} & \texttt{MLR\_startup\_costs} \\ 
 \hline
 Network & \texttt{ptdf\_power\_flow} & \texttt{ptdf\_power\_flow} \\ 
\end{tabular}
\caption{\label{ta:formulation}Formulation used in UC and ED models.}
\end{table}

% As our first example, we simulate RTS-GMLC grid by running UC with day-ahead forecast and ED with both forecast and actual data. For every day in 2020, cost attributions are computed between forecast (baseline) and actual data for a $24$-hour period and then averaged across $24\times365=8760$ hours. In Table \ref{ta:top10}, we report the top $10$ renewable generators with the largest average cost attributions. It is noticeable that the wind generators are ranked at the top. This is due to the capacity of the units and fact that solar generators only produce power during the day. In fact, the capacity for 317-WIND-1, 309-WIND-1 and 319-PV-1 are $799.1$MW, $148.3$MW and $188.2$MW respectively. The renewable generators exist in $5$ sub-areas in the RTS-GMLC grid. We further compute the sum of the generation deficiency (i.e. actual - forecast), the sum of the cost attribution and the average of the shadow price within each sub-area and plot in Figure \ref{fi:act_fcst_attri}.

\begin{figure}[H]
\centering
{\includegraphics[width=1.0\textwidth,height=0.5\textwidth]{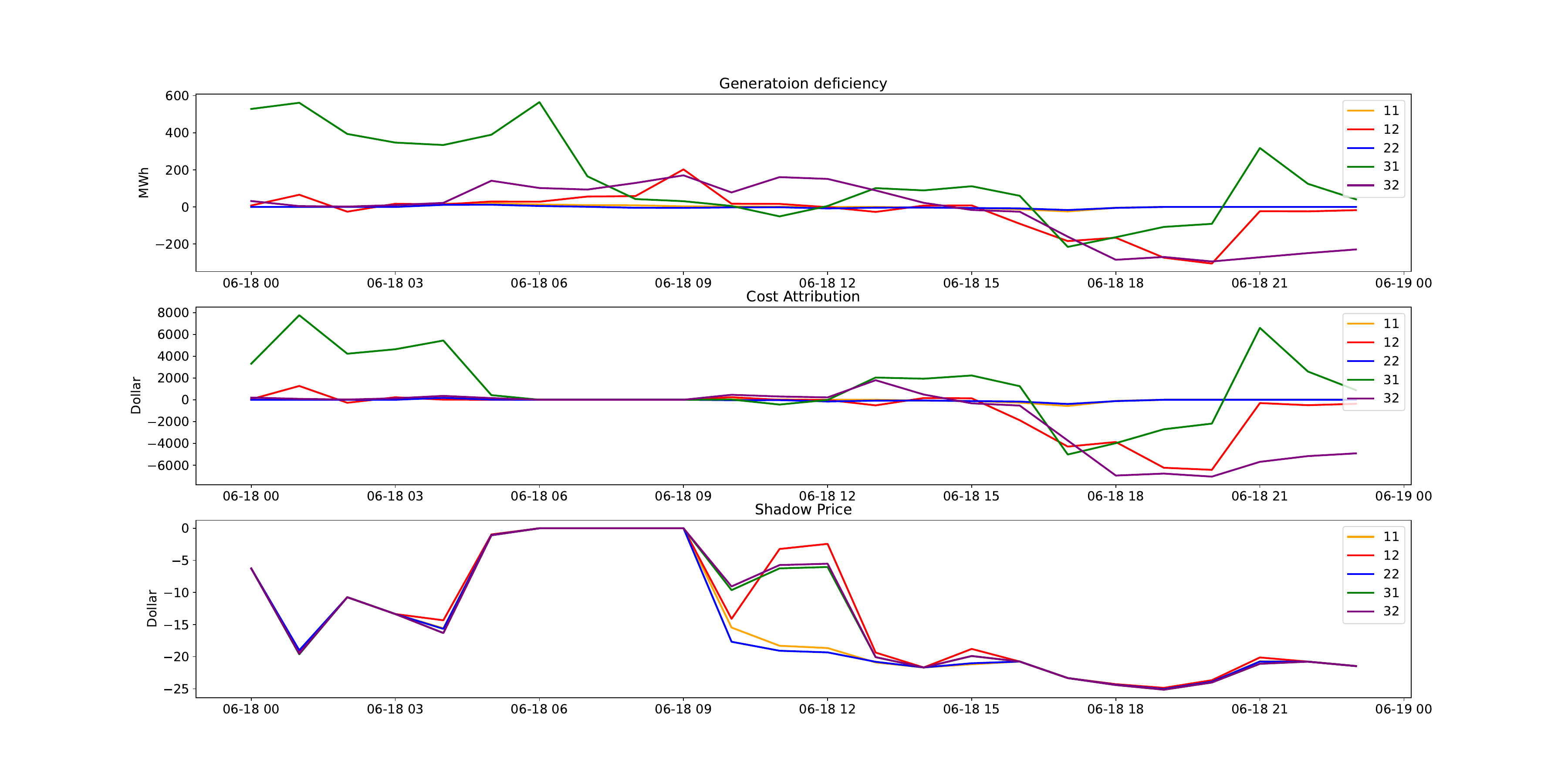}}\quad
\caption{Generation deficiency, cost attributions and shadow prices for sub-areas in RTS-GMLC.}
\label{fi:act_fcst_attri}
\end{figure}

In our initial case study, we simulate the RTS-GMLC grid by executing UC with day-ahead forecasts and ED with both forecasted and actual data. Across each day of 2020, we calculate cost attributions comparing forecasted (baseline) and actual data over 24-hour periods, then average these attributions hourly. Table \ref{ta:top10} presents the top 10 renewable generators with the highest average cost attributions. Notably, wind generators rank highest due to their substantial capacity, followed by solar generators, which operate predominantly during daylight hours. Specifically, the capacities for 317-WIND-1, 309-WIND-1, and 319-PV-1 are $799.1$MW, $148.3$MW, and $188.2$MW, respectively. These renewable generators are distributed across five sub-areas within the RTS-GMLC grid. Additionally, we calculated the sum of generation deficiencies (actual - forecast), the sum of cost attributions, and the average shadow price within each sub-area for 2020-06-18. They are provided in Figure \ref{fi:act_fcst_attri}.

\begin{table}[hbt!]
\begin{tabular}{ c|c|c|c } 
\label{ta:top10}
 Generator & Bus & Attribution & Shadow Price\\ 
  \hline
  317-WIND-1 & 317 & $54625.9$ & $-159.2$ \\
  \hline
  122-WIND-1 & 122  & 	$44742.8$ &	$-151.1$ \\
  \hline
  303-WIND-1 & 303 & $32993.6$ & 	$-148.8$ \\
  \hline
  309-WIND-1 &	309 &	$6318.4$ 	& $-168.7$ \\
  \hline 
  319-PV-1 & 319 & $768.9 $ &	$-161.6$ \\
  \hline 
  215-PV-1 &	215 &	$481.5$ &	$-154.2$\\
  \hline
  113-PV-1 &	113 &	$380.0$ &	$-153.1$\\
  \hline
  313-PV-2 &	313 &	$257.2$ &	$-164.6$\\
  \hline
  312-PV-1 &	312 &	$237.7$ &	$-165.2$ \\
  \hline 
  310-PV-2 &	310 &	$220.2$ & 	$-164.3$\\
\end{tabular}
\caption{\label{ta:top10}Top $10$ renewable generators with the largest average cost attributions.}
\end{table}

Next, we examine the performance of the algorithm using Monte Carlo scenarios. The load demand and renewable production scenarios used for computing the cost attributions are generated by the open source Python package PGscen \cite{pgscen}. For each day in 2020, we first generated $K=1000$ Monte Carlo scenarios for load demand and renewable (wind and solar) production. With the risk level $\alpha=0.05$, the algorithm is then applied to obtain cost attributions for the worst $\alpha K=50$ scenarios.

All the experiments are conducted on a 80-core cluster equipped with two Intel(R) Xeon(R) Platinum 8380 CPU @ 2.30GHz and 512GB memory. For each scenario, the cost attribution procedure is run on a single core with no parallelism. To examine the performance of the cost attribution method, we consider the following relative efficiency gap:

\begin{equation*}
    \text{relative efficiency gap} = \frac{\max_{1\leq\tau\leq24}\left|C_{act}^{\tau} - C_{fcst}^{\tau} - \sum_{m=1}^M C_{init}^{m,\tau} - \sum_{l=1}^L C_{load}^{l,\tau} - \sum_{n=1}^N C_{renew}^{n,\tau}\right|}{\max_{1\leq\tau\leq24}\left|C_{act}^{\tau}\right|}.
\end{equation*}

We report in Table \ref{ta:alloc_summary} the average number of actual quadrature nodes, run time (in seconds) and the relative efficiency gap over $50\times364=18200$ scenarios from 2020-01-02 to 2020-12-30. 

\begin{table}[H]
     \begin{tabular}{c|c|c|c}
      & number of nodes & run time (s) & relative efficiency gap \\
      \hline
      max & $124$ & $1307.9$ & $5.3\%$ \\
      \hline
      median & $23$ & $239.0$ & $0.2\%$ \\
      \hline
      mean & $26.2$ & $275.8$ & $0.6\%$\\
      \hline
      std & $18.0$ & $195.8$ & $0.8\%$ \\
    \end{tabular}
    \caption{\label{ta:alloc_summary}Cost attribution performance.}
\end{table}

\begin{figure}[H]
\centering
{\includegraphics[width=1.0\textwidth,height=.3\textwidth]{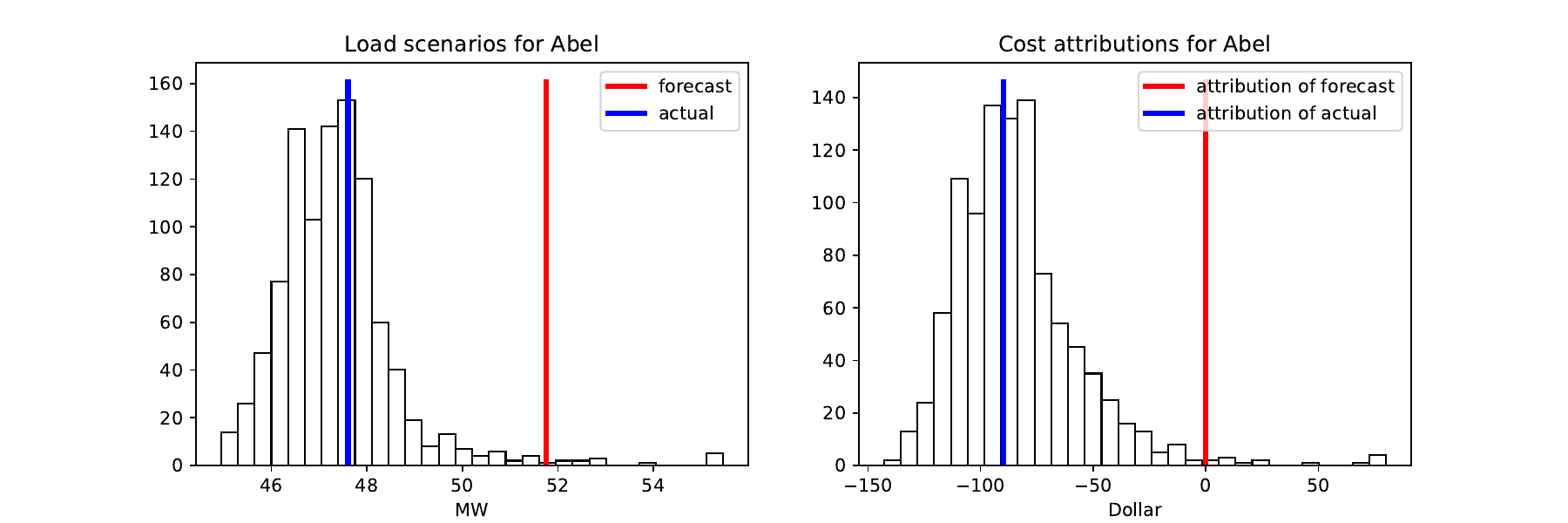}}\quad
\caption{Scenarios and cost attributions for Abel at 2020-07-08 04:00.}
\label{fi:hist_load}
\end{figure}

\begin{figure}[H]
\centering
{\includegraphics[width=1.0\textwidth,height=.3\textwidth]{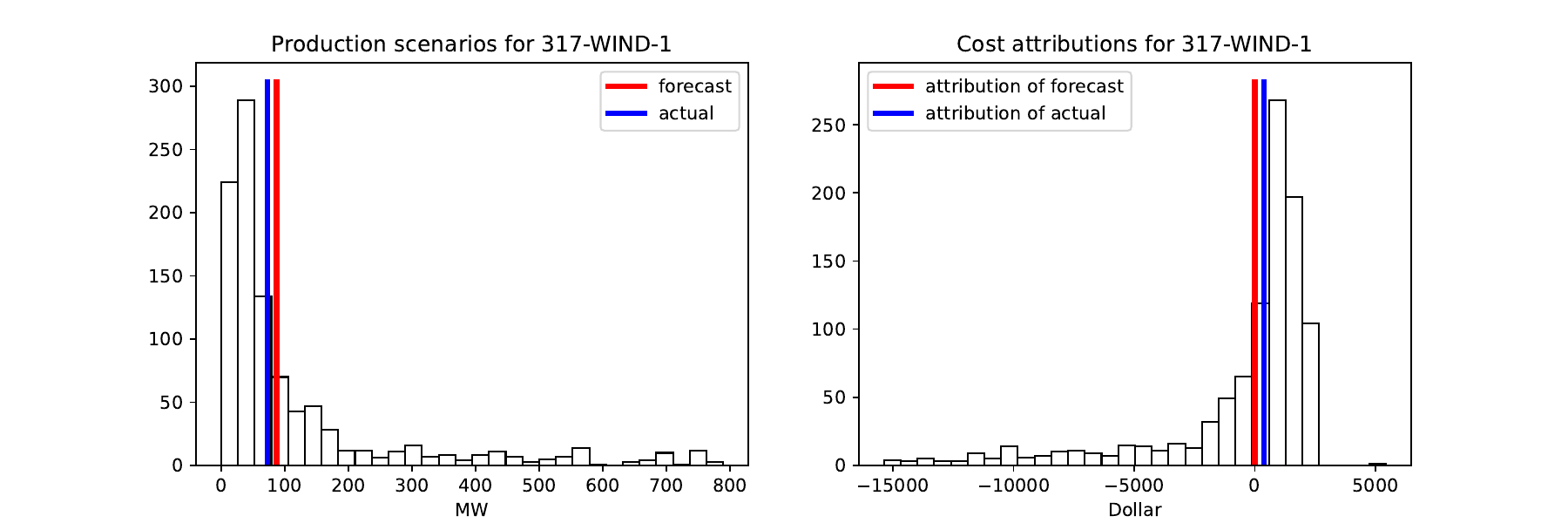}}\quad
\caption{Scenarios and cost attributions for 317-WIND-1 at 2020-07-08 17:00.}
\label{fi:hist_wind}
\end{figure}

\begin{figure}[hbt!]
\centering
{\includegraphics[width=1.0\textwidth,height=.3\textwidth]{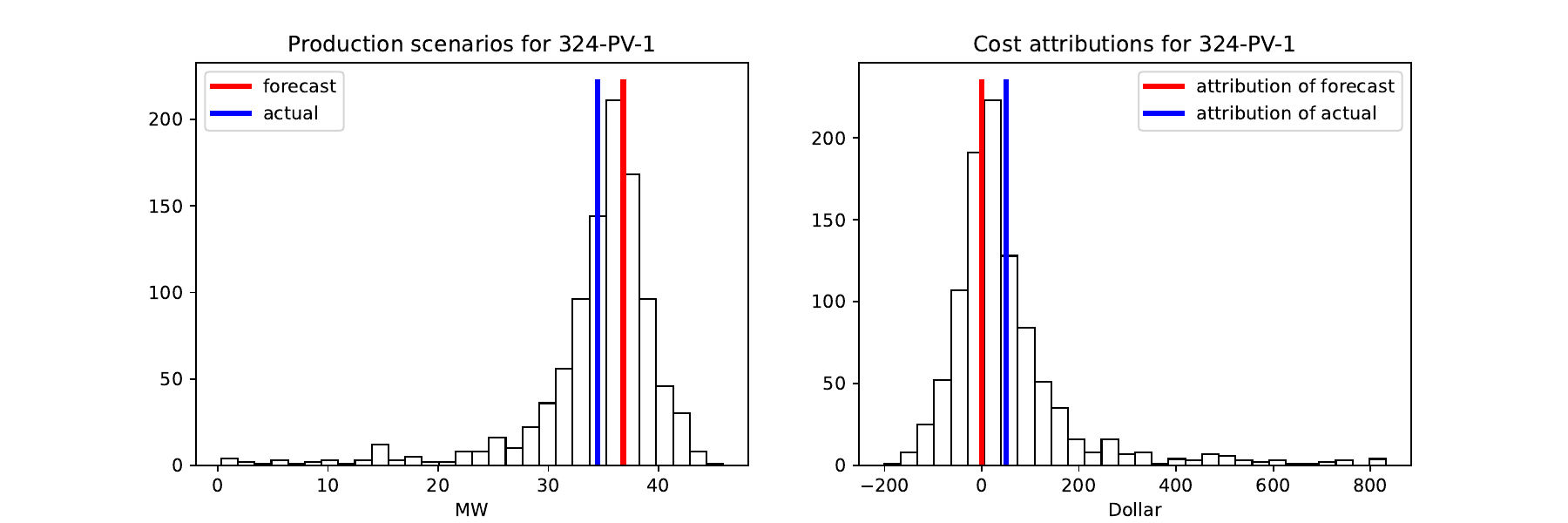}}\quad
\caption{Scenarios and cost attributions for 324-PV-1 at 2020-07-08 15:00.}
\label{fi:hist_solar}
\end{figure}

For the purpose of illustration, we present numerical results for the test day 2020-07-08. We performed the cost attribution computations for $1000$ Monte Carlo scenarios and 
we discuss the results for the load bus Abel, the wind farm 317-WIND-1 and the solar plant 324-PV-1. In each case, we
show in Figures \ref{fi:hist_load}, \ref{fi:hist_wind} and \ref{fi:hist_solar} the histograms of the values of the scenarios (left) and their cost attributions (right) computed scenario by scenario.
Figure \ref{fi:hist_load} shows that on that day, the load forecast significantly over estimated the actual load, and the Monte Carlo generation engine (which does not have a crystal ball and could not use the value of the actual load to generate the scenarios) used a heavy tail distribution to produce reasonable values despite the forecast overshoot.
The two histograms follow similar patterns. This should be expected  due to the positive correlation between load demand and system costs. 

\begin{figure}[H]
\centering
{\includegraphics[width=1.0\textwidth,height=.3\textwidth]{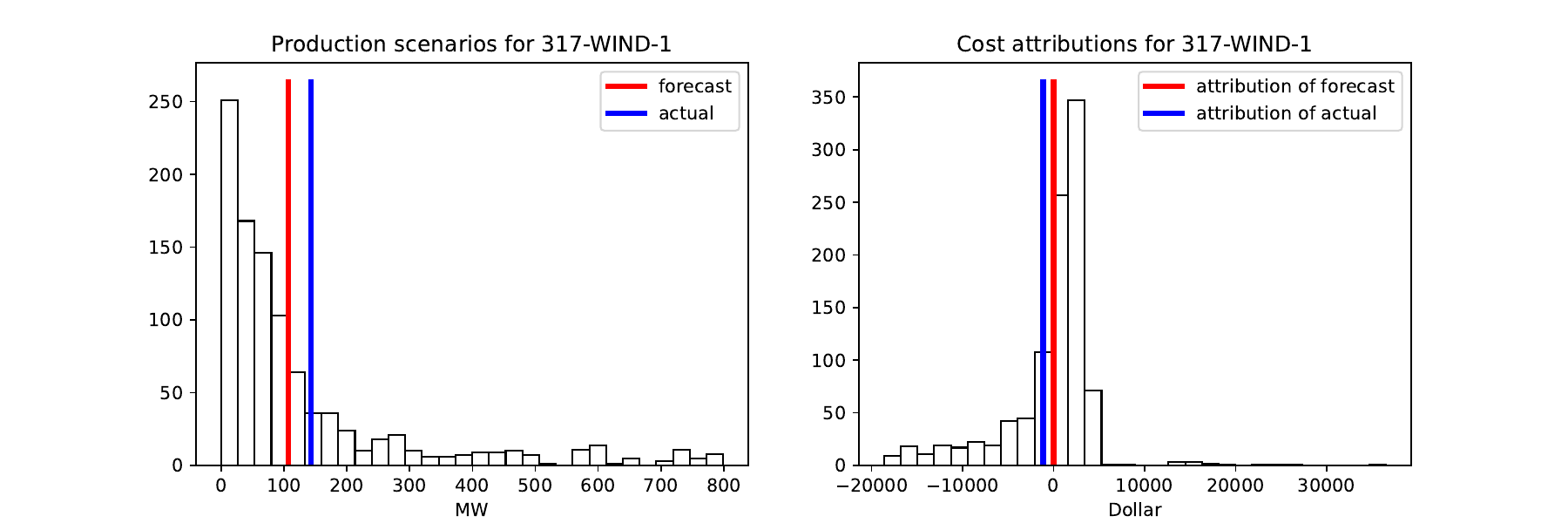}}\quad
\caption{Scenarios and cost attributions for 317-WIND-1 at 2020-07-08 19:00.}
\label{fi:hist_wind_loadshedding}
\end{figure}

Figure \ref{fi:hist_wind} shows that on that day, the wind forecast at this location was pretty accurate. Still, the patterns of the histograms seem to be mirror images of each other. This is consistent  with the negative correlation between power production and system costs.
The same remarks apply to Figure \ref{fi:hist_solar}.  On the other hand, it's worth noting that when a loss of load occurs during a specific hour within a scenario, the attribution of the associated cost is likely to be significantly amplified. This is due to the substantial penalty incurred for failing to meet the load demand. As illustrated in Figure \ref{fi:hist_wind_loadshedding}, we can observe an example of this phenomenon where certain scenarios with low production are attributed with notably high costs. This emphasizes the critical importance of addressing load shedding and its cost implications in our analysis.

%%%%%%%%%%%%%%%%%%%%%%%%%%%%
%%%%%%%%%%%%%%%%%%%%%%%%%%%%%
\section{\textbf{Risk-Averse Unit Commitment}}
\label{se:risk_averse_uc}

\subsection{Scenario-Based Risk Allocation}

In this section, we demonstrate an application of the cost attribution algorithm within a simulation framework tailored for a risk-averse Unit Commitment (UC) approach. First, we introduce a risk allocation algorithm built upon the aforementioned cost attribution method. The primary objective of this algorithm is to assess the impact of uncertainties on daily operational costs of the system. This is achieved by incorporating the variability in both load and renewable generation across numerous Monte Carlo scenarios. Initially, we solve the UC and ED problems using day-ahead forecasts. Subsequently, we use these scenarios as proxies for actual load and renewable production data and solve ED once again. The scenarios employed in this study are generated utilizing the model described in \cite{carmona2022joint}, as implemented in \cite{pgscen}. This model captures a high-dimensional spatial-temporal correlation structure and encompasses extreme and risky events.

To be more specific, let us assume that $K$ Monte Carlo scenarios have been generated and a risk level $\alpha\in\left(0,1\right)$ has been chosen. The risk allocation algorithm comprises the following steps:

\begin{enumerate}
    \item The UC and ED are solved based on the day-ahead forecasted data.
    \item For each scenario, its values are substituted for the actual data  and we solve the $24$ hourly ED models to obtain hourly costs.
    \item We select the worst $\alpha K$ scenarios in terms of their total cost. We denote the set of the worst scenarios by $\mathcal{S}$.
    \item For each scenario in $\mathcal{S}$, we calculate the cost attribution between the baseline model, which uses forecasted data, and the target model, which uses scenario data. This enables us to assign costs to all initial states of dispatchable generators, load demands, and renewable generators for every hour.  
    \item \label{step:average cost} For a given asset (load demand or renewable generator) at any hour, we compute its risk score $R$ as the average of its cost attributions across the scenarios in $\mathcal{S}$. For example, the risk score for the renewable generator $n$ at hour $\tau$ is $R=\frac{1}{\alpha K}\sum_{j\in\mathcal{S}} C_{renew}^{n,\tau, j}$, where $C_{renew}^{n,\tau, j}$ is the cost attribution in the $j$-th scenarios.
\end{enumerate}

%%%%%%%%%%%%%%%%%%%%%%
\subsection{Reliability Adjustment}
As previously explained, our proposed method assigns a reliability risk score to each stochastic input during every hour of the day. While these stochastic inputs include various factors such as loads and outputs from renewable generation assets, our attention in this study will be concentrated on the renewable generation assets. To mitigate the risk associated with the uncertainty of these renewable generators, our proposal is to modify their generation capacities whenever elevated risk levels are identified.

\begin{figure}[hbt!]
\centering
{\includegraphics[width=.45\textwidth,height=.3\textwidth]{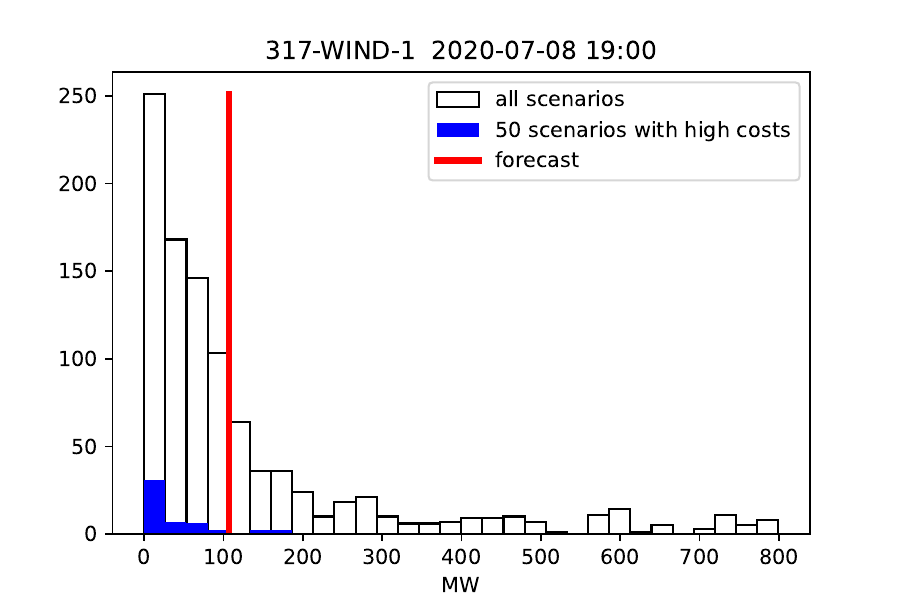}}{\includegraphics[width=.45\textwidth,height=.3\textwidth]{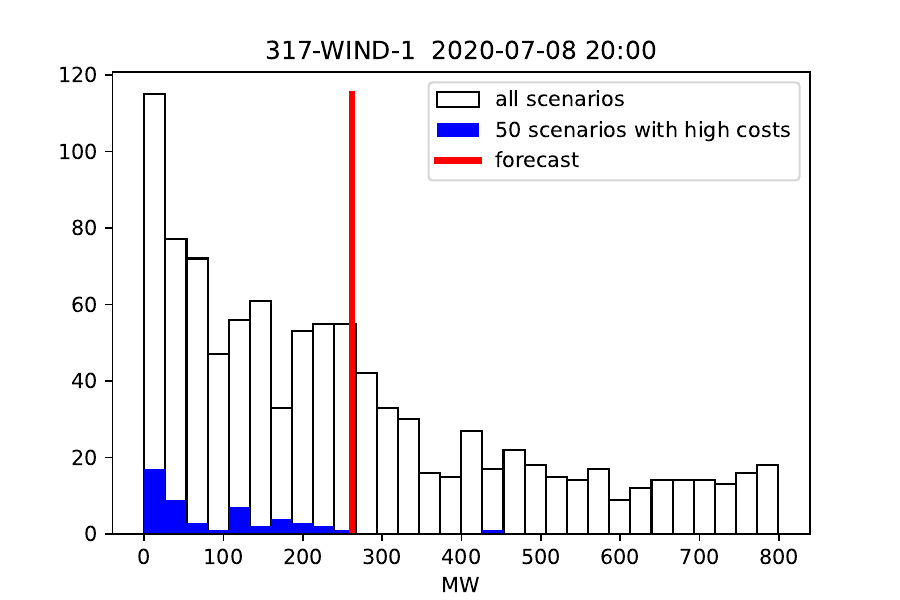}}
\caption{Scenarios and the worst case scenarios for two wind farms on 2020-04-26.}
\label{fi:renew_adjust_hist}
\end{figure}

Specifically, given a risk level $\alpha\in(0,1)$, we recall that the risk score $C_{renew}^{n,\tau}$ for the renewable generator $n$ at hour $\tau$ is computed as the average costs attributions over a set of worst case scenarios. These costs are frequently linked to the under-production of the renewable generators relative to the forecasts. Consequently, our initial step involves the computation of a "per MWh" risk score $R_{\tau}^n$ through

\begin{equation*}
    R_{\tau}^n = \frac{C_{renew}^{n,\tau}}{q_{\tau}^{n, fcst}-\frac{1}{K}\sum_{j\in\mathcal{S}}q_{\tau}^{n,j}}
\end{equation*}
where $q_{\tau}^{n,j}$ is the generation capacity of the renewable generator $n$ at hour $\tau$ in the $j$-th scenario. Following that, we introduce two positive parameters, denoted as $\underline{R}$ and $\overline{R}$. Broadly speaking, these parameters are meant to represent the per MWh cost of the grid under usual conditions and in a situation where the generation needs to be completely disregarded. Subsequently, a generation capacity adjustment percentage can be calculated as follows:

\begin{equation*}
    r_{\tau}^n = \max\left(0, \min\left(1, \frac{R_{\tau}^n-\underline{R}}{\overline{R}}\right)\right).
\end{equation*}
where $\min$ and $\max$ functions guarantees that $r_{\tau}^n\in\left[0, 1\right]$. Finally, the adjusted renewable generation capacity is determined by the equation:  
\begin{equation*}
    q_{\tau}^{n,adj} = q_{\tau}^{n,fcst} - r_{\tau}^n\left(q_{\tau}^{n,fcst}-q_{\tau}^{n,min}\right).
\end{equation*}
Here $q_{\tau}^{n,min}$ represents the minimum (\emph{guaranteed}) generation capacity of renewable generator $n$ at hour $\tau$ across all scenarios. By design, when $r_{\tau}^n=0$, indicating no need for adjustment, the equation simplifies to $q_{\tau}^{n,adj} = q_{\tau}^{n,fcst}$. Conversely, if $r_{\tau}^n=1$, implying a necessity for maximum adjustment, the generation capacity is reduced to the minimum (\emph{guaranteed}) level: $q_{\tau}^{n,adj} = q_{\tau}^{n,min}$.
As an illustrative example, we demonstrate the generation capacity adjustments for 317-WIND-1 at 2020-07-08 during the hours 19:00 and 20:00. For this example, we choose $\underline{R}=20$ and $\overline{R}=200, 300$, and $500$. The histograms of generation scenarios are presented in Figure \ref{fi:renew_adjust_hist}, and the relevant quantities for computing the adjustments are detailed in Table \ref{ta:renew_adjust}.

\begin{table}
  \begin{tabular}{|c|c|c|c|c|c|c|c|c|c|}
  \hline
    
    \multicolumn{1}{|c|}{\multirow{2}{*}{Hour}} &  
    \multicolumn{1}{|c|}{\multirow{2}{*}{$q_{\tau}^{n,fcst}$}} & 
    \multicolumn{1}{|c|}{\multirow{2}{*}{$C_{renew}^{n,\tau}$}} & \multicolumn{1}{|c|}{\multirow{2}{*}{$R_{\tau}^{n}$}} &
      \multicolumn{2}{|c|}{$\overline{R}=200$} &
      \multicolumn{2}{|c|}{$\overline{R}=300$} &
      \multicolumn{2}{|c|}{$\overline{R}=500$} \\ \cline{5-10}
      \multicolumn{1}{|c}{} & \multicolumn{1}{|c}{} & \multicolumn{1}{|c}{} & \multicolumn{1}{|c}{} &
      \multicolumn{1}{|c}{$r_{\tau}^n$} &
      \multicolumn{1}{|c}{$q_{\tau}^{n,adj}$} &
      \multicolumn{1}{|c}{$r_{\tau}^n$} &
      \multicolumn{1}{|c|}{$q_{\tau}^{n,adj}$} &
      \multicolumn{1}{|c}{$r_{\tau}^n$} &
      \multicolumn{1}{|c|}{$q_{\tau}^{n,adj}$}\\
    \hline
    \hline
    \multirow{1}{*}{$19$} & $107.0$ & 
   $15630.2$ & $210.8$	 & $95.4\%$ & $4.9$ & $63.6\%$ & $38.9$ & $ 38.2\%$ & $66.2$\\
    \hline
    \multirow{1}{*}{$20$} & $262.0$ & 
    $6163.0$ & $35.0$ & $7.5\%$ & $242.3$ & $5.0\%$ 
    & $248.9$ & $3.0\%$ & $254.1$\\
    \hline
  \end{tabular}
  \caption{\label{ta:evscore}Generation capacity adjustments for 317-WIND-1 at hours $19$ and $20$ of 2020-07-08.}
  \label{ta:renew_adjust}
\end{table}

%%%%%%%%%%%%%%%%%%%%%%%%%
\subsection{\textbf{Simulation results}}

\begin{figure}[H]
\centering
{\includegraphics[width=1.0\textwidth,height=.6\textwidth]{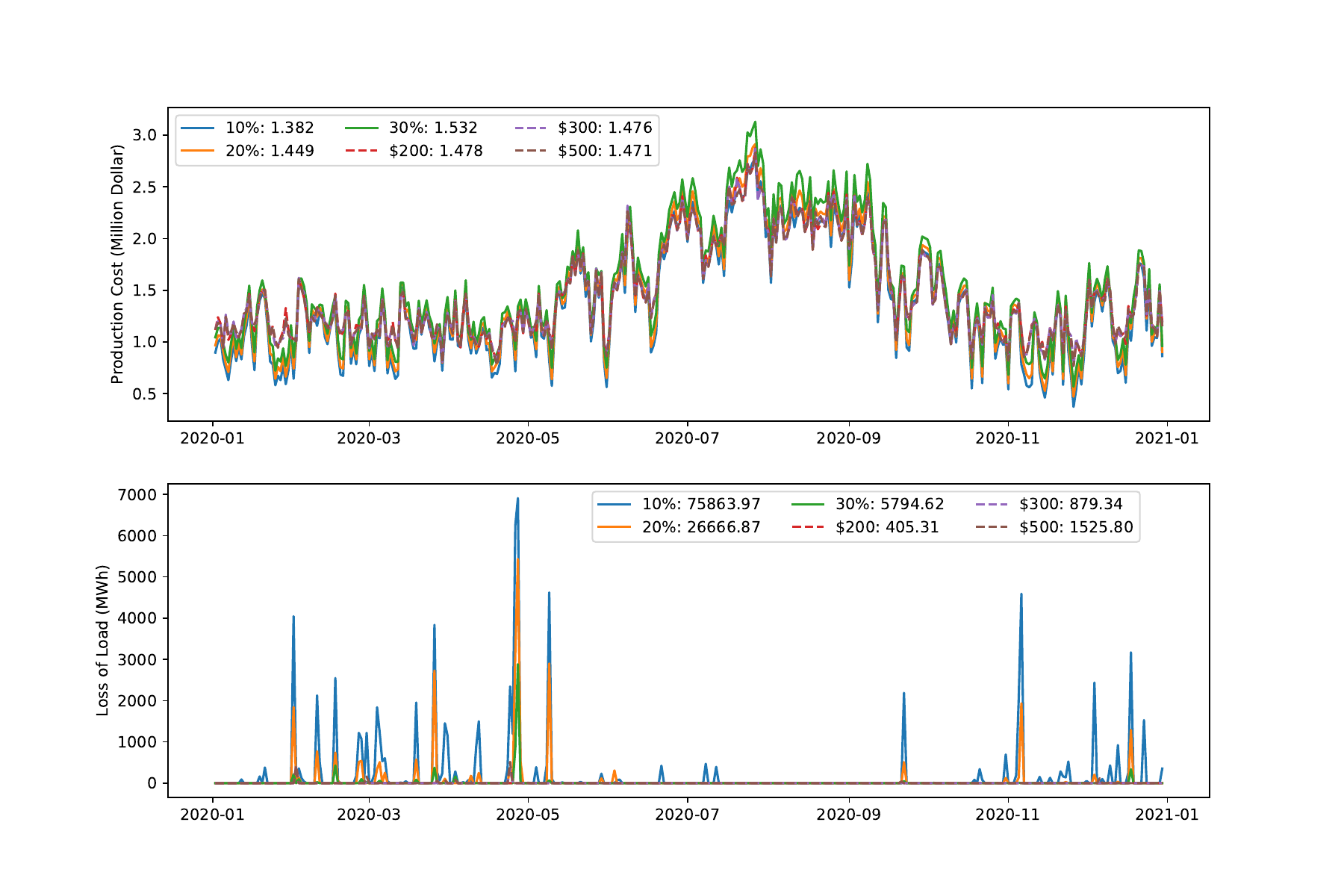}}\quad
\caption{Production costs and loss of load in RTS-GMLC system using different UC from 2020-01-02 to 2020-12-30.}
\label{fi:renew_adjust_uc_yearly}
\end{figure}

To examine the impact of the proposed risk-averse UC approach, we run simulations for the RTS-GMLC grid from 2020-01-02 to 2020-12-30. For the purpose of comparison, in addition to our approach, we include three benchmark approaches, in which the standard alternating UC and ED are performed with a system spinning reserve equal to $10\%$, $20\%$ and $30\%$ of the system forecasted load demand. In the benchmark approaches, the load mismatch (over-generation and loss of load) and reserve shortfall penalties in the objective function of the UC model are selected to be $10000$\$/MWh and $1000$\$/MWh. In contrast, our approach adopted a reserve factor of $5\%$ and used the same load mismatch and reserve shortfall penalties as the benchmark methods. For generation capacity adjustments, $\underline{R}$ was set to $20\$$, and the value of $\overline{R}$ was varied by $200\$$, $300\$$, and $500\$$. In Figure \ref{fi:renew_adjust_uc_yearly}, we show the daily system operational costs and the load shedding patterns for the entire year using different UC approaches. Additionally, we report the average daily operational costs and total yearly load shedding in the same plot. We observe that the proposed risk-averse UC approach is capable of largely reducing load shedding (which is the largest risk to the grid operation) under all choices of $\overline{R}$. Specifically, selecting the least aggressive parameter $\overline{R}=500$ results in a $1525$ MWh loss of load, roughly a quarter of the load shedding experienced with a $30\%$ reserve factor. The production costs using the risk-averse UC approach are approximately \$$1.47$ million, which is much lower than the \$$1.53$ million with a $30\%$ reserve factor and only slightly higher than with a $20\%$ reserve factor.

\begin{figure}[H]
\centering
{\includegraphics[width=1.0\textwidth,height=.6\textwidth]{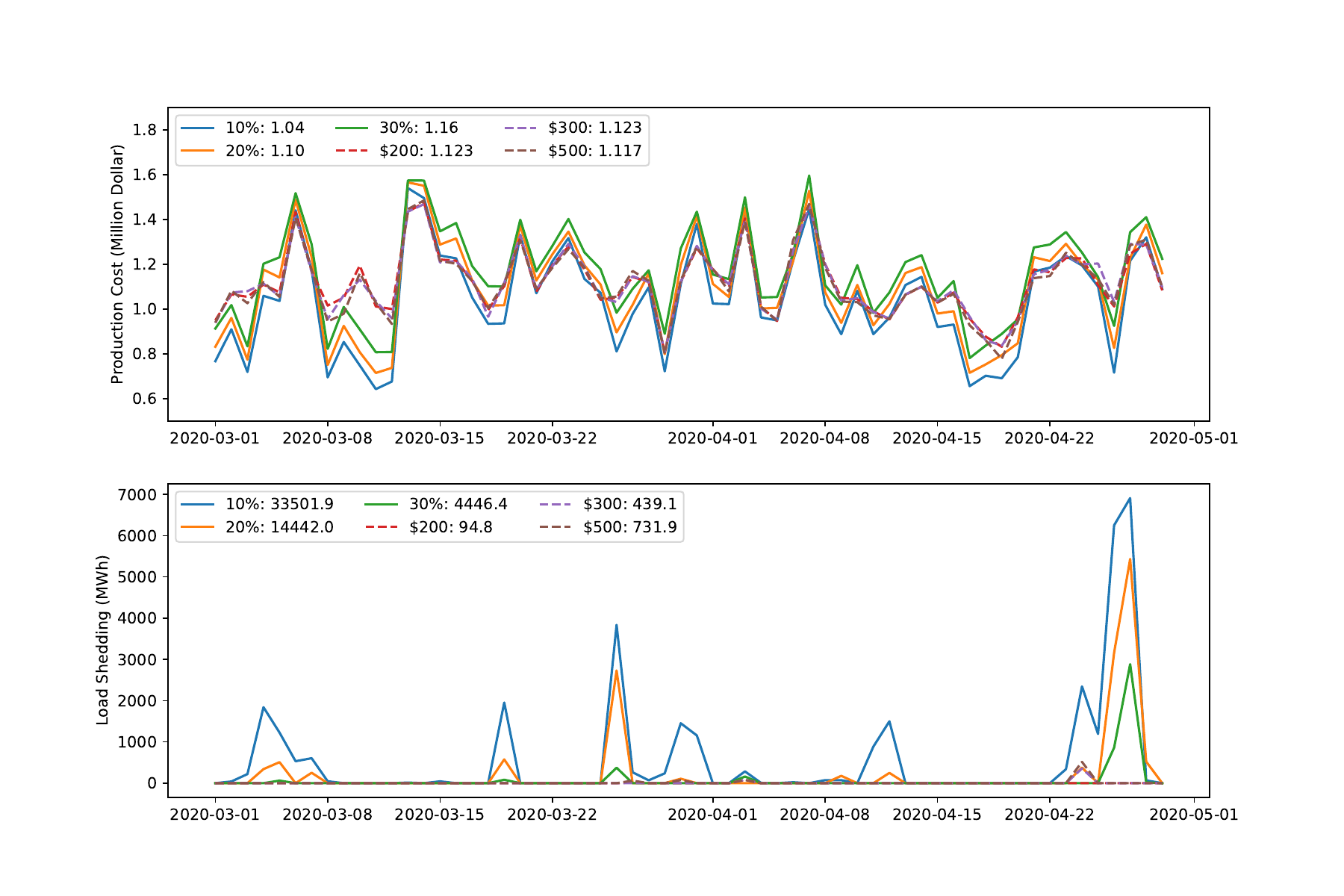}}\quad
\caption{Production costs and loss of load in RTS-GMLC system using different UC for 2020-03-01 to 2020-04-30.}
\label{fi:renew_adjust_uc_monthly}
\end{figure}

In Figure \ref{fi:renew_adjust_uc_monthly}, we create the same plot from 2020-01-01 to 2020-04-30, a shorter period in which activities like load shedding happen more frequently. For 2020-03-26, 2020-04-26 and 2020-04-27, it is clear that the risk-averse UC hedges the risk against the worst case scenarios by depending less on the renewable generations (and thus higher production cost) to avoid possible load shedding due to uncertainties in the day-ahead forecasts. The loss of load for all choices of $\overline{R}$ is less than $1000$ MWh, whereas reserve factors of $10\%$, $20\%$, and $30\%$ result in unmet load demands ranging from $4000$ to $34000$ MWh. The production costs using the risk-averse UC approach, however, are comparable to those incurred with a $20\%$ reserve factor.

\subsection{Load Loss Abatement}
\begin{figure}[H]
\centering
\subfloat[][]{\includegraphics[width=0.85\textwidth,height=.6\textwidth]{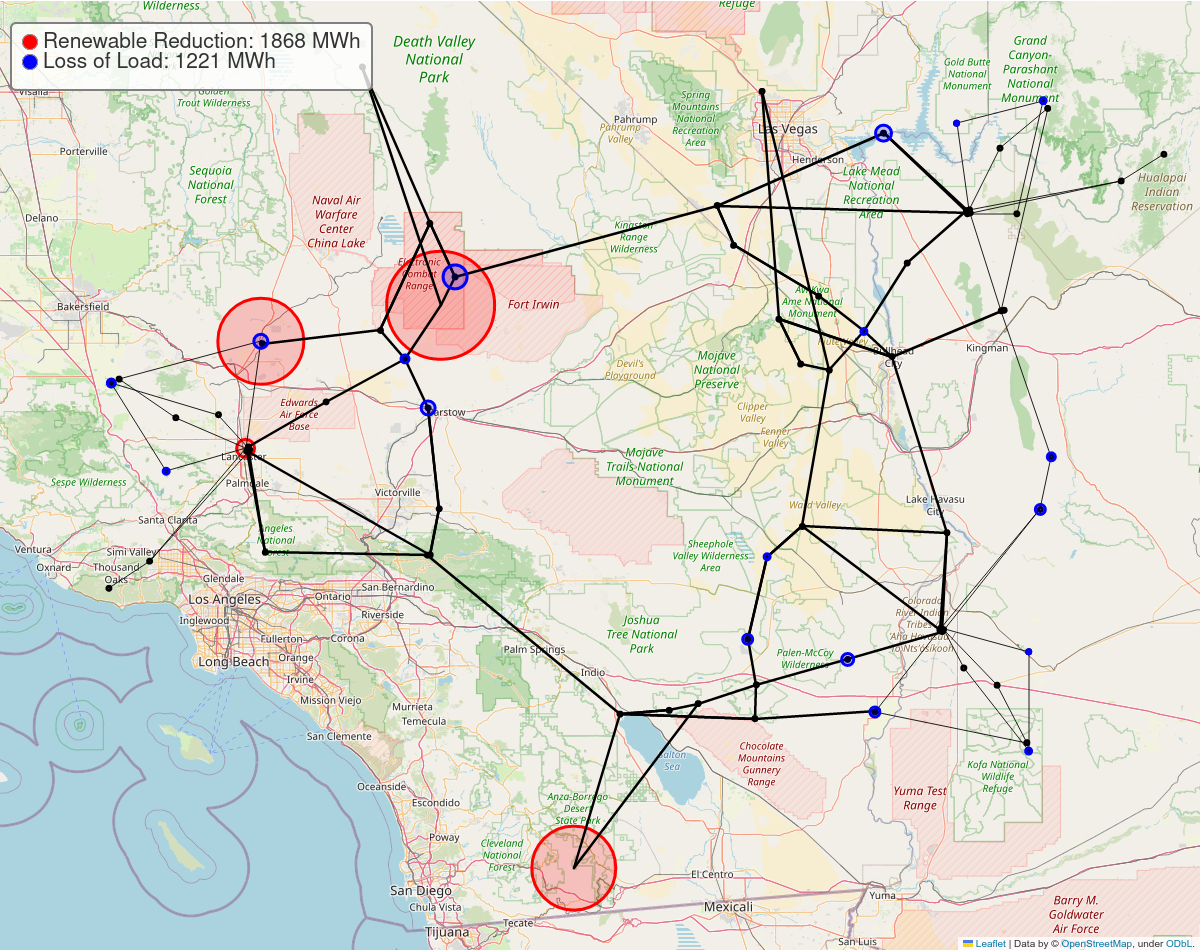}}\\
\subfloat[][]{\includegraphics[width=0.85\textwidth,height=.6\textwidth]{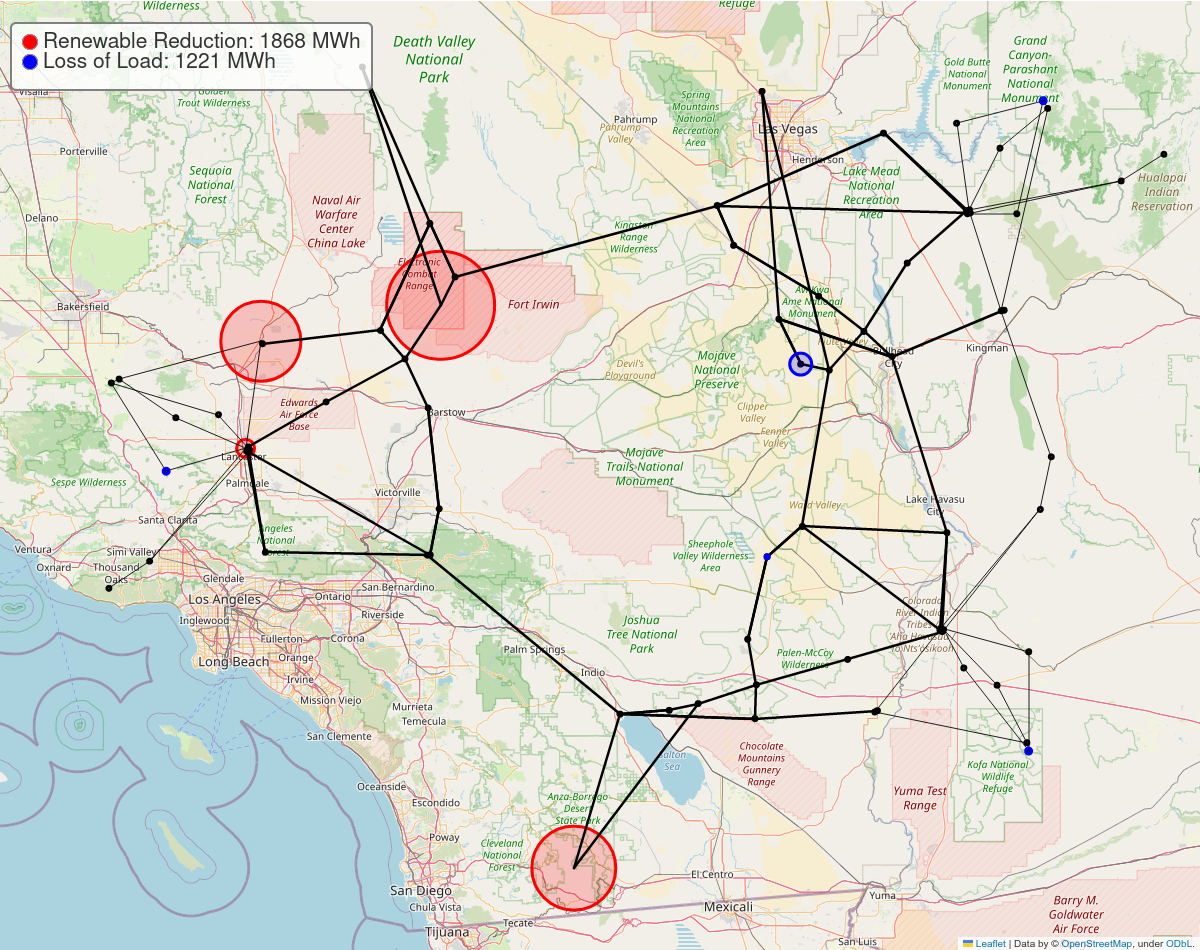}}
\caption{Comparison of risk-neutral UC using $5\%$ (A) and $30\%$ (B) reserve requirements with risk-averse at 2020-04-26 21:00:00.}
\label{fi:loadshedding_map}
\end{figure}

To further explore how the renewable generation adjustment assists in mitigating risk, we present a more comprehensive comparison with the use of a fixed percentage of the load demand spinning reserves factor (referred to as the risk-neutral UC) on 2020-04-26. During the $24$-hour period, when operating the ED with the risk-averse UC, there were no instances of load loss. However, running it with the risk-neutral UC using $5\%$ and $30\%$ reserve factors resulted in a total load loss of approximately $860$ and $7903$ MWh, respectively. This loss primarily occurred due to an overestimation of wind power generation in the day-ahead forecasts. The day-ahead predictions estimated wind power plants to generate as much as 37,046 MWh throughout the day, while the actual real-time production was only $12590$ MWh. Notably, at 21:00, there was an over-forecast of wind power production by about $2019$ MWh. In Figure \ref{fi:loadshedding_map}, we illustrate the RTS-GMLC grid on a geographical map, indicating buses with load shedding (depicted as blue circles) in simulations using the risk-neutral UC. Additionally, we calculate the difference in renewable power dispatch in simulations using the risk-averse UC, revealing a significant reduction in renewable power dispatch due to the associated capacity adjustments, shown as differences represented by red circles on the map.

%%%%%%%%%%%%%%%%%%%%%%%%%%%%

%%%%%%%%%%%%%%%%%%%%%%%%%%%%%%%%%%%%%%%%%%%%%%
\section{\textbf{Conclusion}}
\label{se:conclusion}

In conclusion, this paper offers a new perspective on addressing uncertainty and associated risks in power system management, driven by the increasing integration of renewable energy sources. We propose a novel algorithm that leverages the discrepancies between forecasted and actual values to quantify the financial risks associated with uncertainty. Inspired by Integrated Gradients (IG), a renowned method in machine learning, our algorithm attributes the contributions of stochastic components to differences in system costs. By treating the power production cost model as a function that takes inputs of random quantities like load demand and renewable generation, we provide a practical approach for understanding the impact of variability on grid operational costs.

Moreover, we demonstrate the utility of our approach in a risk-averse unit commitment framework, where adjustments to renewable generator capacities are made based on Monte Carlo scenarios. This framework mitigates system risk by reducing reliance on highly uncertain renewable productions, thereby enhancing grid reliability. On the other hand, the overall operational costs will not increase significantly under this framework (see Figures \ref{fi:renew_adjust_uc_yearly} and \ref{fi:renew_adjust_uc_monthly}) as it identifies critical times and locations where highly uncertain renewable productions occur.

The results of our simulations on the RTS-GMLC system serve as a compelling validation of the effectiveness of our algorithms. Furthermore, they highlight the superior performance of the risk-averse UC framework, particularly in its ability to significantly reduce the occurrence of load shedding. As part of our future work, we plan to further refine and extend our approach, conducting tests on larger and more realistic grid systems to verify its performance in more complex and practical scenarios.

%%%%%%%%%%%%%%%%%%%%%%%%%%%%%%%%%%%%%%%%%%%%%%

%%%%%%%%%%%%%%%%%%%%%%%%%%%%%%%%%%%%%%%%%%%%%
\section{\textbf{Declaration of Competing Interest}}
%%%%%%%%%%%%%%%%%%%%%%%%%%%%%%%%%%%%%%%%%%%%%
The authors declare that they have no known competing financial interests or personal relationships that could have appeared to influence the work reported in this paper.

%%%%%%%%%%%%%%%%%%%%%%%%%%%%%%%%%%%%%%%%%%%%%
\section{\textbf{Acknowledgments}}
%%%%%%%%%%%%%%%%%%%%%%%%%%%%%%%%%%%%%%%%%%%%%
The authors were partially supported by ARPA-E grants DE-AR0001289 and DE-AR0001390 under the PERFORM program of the US Department of Energy. We would like to thank Mike Ludkovski (University of California at Santa Barbara) and Glen Swindle (Scoville Risk Partners) for enlightening conversations on the content of the paper; and Michal Grzadkowski (Princeton University) for software development and implementation.

%%%%%%%%%%%%%%%%%%%%%%%%%%%%%%%%%%%%%%%%%%%%%
\bibliographystyle{plain}
% \small
%\bibliography{PERFORM.bib}

\end{document}